\documentclass[final]{dmtcs-episciences}

\usepackage{amsfonts,amsmath,amssymb,amsthm}
\usepackage{cleveref}
\usepackage{mathtools}
\newtheorem{theorem}{Theorem}
\newtheorem{lemma}[theorem]{Lemma}
\newtheorem{example}[theorem]{Example}

\newcommand*{\PF}{\mathrm{PF}}
\newcommand*{\PP}{\mathrm{PP}}
\newcommand{\seqnum}[1]{\href{https://oeis.org/#1}{#1}}

\usepackage[utf8]{inputenc}
\usepackage{subfigure}

\usepackage[round]{natbib}

\author[Harris, Mart\'inez Mori, and Wilson]{Pamela E. Harris\affiliationmark{1}
  \and J. Carlos Mart\'inez Mori\affiliationmark{2}
  \and Alexander N. Wilson\affiliationmark{3}}
\title[A Pollak Proof for the Number of Weakly Increasing Parking Functions]{A Pollak Proof for the Number of Weakly Increasing Parking Functions}
\affiliation{
  Department of Mathematical Sciences, University of Wisconsin-Milwaukee, Milwaukee, WI, USA\\
  Department of Mathematical and Statistical Sciences, University of Colorado Denver, Denver, CO, USA\\
  Department of Mathematics and Statistics, York University, Toronto, ON, Canada}
\begin{document}
\publicationdata{vol. 28:1, Permutation Patterns 2025}{2026}{2}{10.46298/dmtcs.17006}{2025-12-01; 2025-12-01; 2026-04-03}{2026-04-07}

\maketitle
\begin{abstract}
We develop a circular-street argument, in the style of Pollak, to obtain a new proof that there are $C_n = \frac{1}{n+1}\binom{2n}{n}$ weakly increasing parking functions of length $n \geq 1$, where $C_n$ is the $n$th Catalan number.
\end{abstract}

\section{Introduction}
\label{sec: introduction}
Consider a one-way street with $n \in \mathbb{N} \coloneqq \{1, 2, 3, \ldots \}$ numbered parking spots.
A sequence of $n$ cars arrive to park on the street one at a time.
When car $i \in [n] \coloneqq \{1, 2, \ldots, n\}$ arrives, it immediately drives to its preferred spot $a_i \in [n]$.
If the spot is empty, then the car parks there.
Otherwise, the car continues driving down the street and parks in the first unoccupied spot it finds after its preference, if any.
If the car does not find any unoccupied spots, then it exits the street unable to park.
Let the tuple $\alpha = (a_1, a_2, \ldots, a_n) \in [n]^n$ summarize the cars' parking preferences; we refer to it as a \emph{parking preference tuple}.
If all cars are able to park given the parking preferences in $\alpha$, then $\alpha$ is said to be a \emph{parking function} of length $n$.
Let $\PF_n \subseteq [n]^n$ denote the set of parking functions of length $n \geq 1$.

Parking functions were first studied by \citet{pyke1959supremum} and by \citet{konheim1966occupancy}; they independently showed that $|\PF_n| = (n+1)^{n-1}$ for all $n \geq 1$ \cite[\seqnum{A000169}]{OEIS}.
\citet{riordan1969ballots} credits Pollak with a different and particularly elegant proof of this result based on a circular-street argument; refer to \citep{martinezmori2025what} for an expository description of Pollak's technique. 

Parking functions are also closely related to Catalan combinatorics; refer to \citet{stanley2015catalan} for a comprehensive survey of this area.
In particular, a parking function $\alpha = (a_1, a_2, \ldots, a_n) \in \PF_n$ is said to be weakly increasing if $a_1 \leq a_2 \leq \cdots \leq a_n$.
Let $\PF_n^\uparrow \subseteq \PF_n$ denote the set of weakly increasing parking functions of length $n \geq 1$.
Then, it is well-known that $|\PF_n^\uparrow| = C_n = \frac{1}{n+1}\binom{2n}{n}$, where $C_n$ is the $n$th Catalan number \cite[\seqnum{A000108}]{OEIS}.
This result can be established in various ways, most notably directly through the recursive relation $C_{n} = \sum_{i=1}^n C_{i} C_{n-i}$ for $n \geq 1$ with $C_0 = 1$ or bijectively through another known Catalan object (most typically Dyck paths, as in \cite[Section 2.2]{armstrong2016rational}).

It is also well-known that, for any $\alpha = (a_1, a_2, \ldots, a_n) \in [n]^n$ with $n \geq 1$, we have that
\begin{equation}
\label{eq: characterization}
    \alpha \in \PF_n \textrm{ if and only if } |\{i \in [n] : a_i \geq j\}| \leq n - j + 1, \ \textrm{ for all } j \in [n].
\end{equation}
This inequality-based characterization implies that parking functions are invariant under the action of the symmetric group, which permutes the entries of a preference tuple: we leverage this property later in this note.
Lastly, as a notational convention, let $[n]_0 \coloneqq [n] \cup \{0\}$ and $\mathbb{N}_0 \coloneqq \mathbb{N} \cup \{0\}$.

Another way of establishing the count of a set of Catalan objects is through a \emph{cyclic shift} argument.
For example, \citet{raney1960functional} developed this kind of proof to count ballot sequences.
In fact, this is part of a long tradition of proofs that partition combinatorial objects into equivalence classes through a rotation action; the cycle lemma of \citet{dvoretzky1947problem} and \citet{golomb56}'s proof of Fermat's little theorem are two additional examples.
Pollak's proof that $|\PF_n| = (n+1)^{n-1}$ is another example of a cyclic shift argument; one that makes specific use of the combinatorial interpretation of parking functions.
However, to the best of our knowledge, no proof in the style of Pollak has been given for the count of weakly increasing parking functions, i.e., for $|\PF_n^\uparrow| = C_n = \frac{1}{n+1}\binom{2n}{n}$.
In \Cref{theorem: catalan} we provide such a proof, along with an illustration of its main ideas in \Cref{ex: proof} and its accompanying \Cref{fig: proof}.
Although our proof is a cyclic shift proof, we distinguish it from other cyclic shift proofs for the enumeration of Catalan objects in that we only employ self-contained parking-related arguments.
In other words, not only is it based on a cyclic shift argument, but more specifically it is based on a circular-street argument.

Our proof makes use of the \emph{content} of a parking preference tuple, i.e., a tuple that records the number of cars that prefer each spot.
Whereas \citep{riordan1969ballots} proved that the number of contents (and therefore the number of weakly increasing preference tuples) that meet \eqref{eq: characterization} satisfies the recurrence relation of the ballot numbers, we instead use a cyclic shift argument to partition this set into equivalence classes each leading to a unique weakly increasing parking function.
The main difficulty lies in that a direct application of Pollak's preference rotations breaks the weakly increasing condition: we correct for this by mapping back into the set of weakly increasing preferences while preserving the decision property of being a parking function.

\section{A Circular-Street Argument for Weakly Increasing Parking Functions}

Consider a sequence of $n$ cars parking on a circular one-way street with $n + 1$ parking spots, numbered $1, 2, \ldots, n, n + 1$.
Let $\PP_n \coloneqq [n + 1]^n$ be the set of all possible parking preference tuples for the $n$ cars.
Similarly, let
\begin{equation*}
   \PP_n^\uparrow \coloneqq \{ \alpha = (a_1, a_2, \ldots, a_n) \in \PP_n : a_1 \leq a_2 \leq \cdots \leq a_n  \}
\end{equation*}
be the set of all possible \emph{weakly increasing} parking preference tuples for the $n$ cars.
Before we present our main result, we present two standalone lemmas that highlight its main ingredients.
\begin{lemma}
\label{lemma: balls and bars}
For any $n \geq 1$, $|\PP_n^\uparrow| = \binom{2n}{n}$.
\end{lemma}
\begin{proof}
This follows as a direct application of the balls and bars method.
Specifically, consider the placement of $n$ bars and $n$ balls on a straight line from left to right.
In this order, the $i$th ball corresponds to the $i$th car, whose preference is $1$ plus the total number of bars to its left.
Therefore, the choice of $n$ places for the bars from the total of $2n$ different places fully determines the cars' weakly increasing preferences.
\end{proof}

\begin{lemma}
\label{lemma: empty spot}
Let $\alpha = (a_1, a_2, \ldots, a_n) \in \PP_n$ and $\alpha' = (a_1', a_2', \ldots, a_n')$ be any rearrangement of $\alpha$.
Then, under either preference tuple $\alpha$ or $\alpha'$, all cars park and ultimately leave the same, single spot unoccupied.
\end{lemma}
\begin{proof}
Since there are $n + 1$ spots but only $n$ cars, regardless of the preference tuple, all cars park and ultimately leave one spot unoccupied.
Let $i^* \in [n + 1]$ be the spot ultimately left unoccupied if the cars arrive under $\alpha$.
We claim that this is the same spot ultimately left unoccupied if the cars arrive under $\alpha'$.
To show this, consider the bijective map $\sigma: [n + 1] \to [n + 1]$ where
\begin{equation*}
    \sigma(j) 
    \coloneqq
    \begin{cases} 
        j + (n + 1 - i^*), &\textrm{ if } 1 \leq j < i^*, \\
        n + 1,  &\textrm{ if } j=i^*, \textrm{ and} \\
        j - i^*, &\textrm{ if } i^* < j \leq n + 1,
    \end{cases}  
\end{equation*}
for all $j \in [n + 1]$.
Essentially, $\sigma$ relabels the spots so that spot $i^* + 1$ becomes the new first spot, spot $i^* + 2$ becomes the new second spot, and so on.
Let $\sigma(\alpha) = (\sigma(a_1), \sigma(a_2), \ldots, \sigma(a_n))$ and note that, based on the assumption that $i^*$ is the spot ultimately left unoccupied under $\alpha$, $n + 1$ is the spot ultimately left unoccupied under $\sigma(\alpha)$.
It follows that $\sigma(\alpha) \in \PF_n$.
Next, consider $\sigma(\alpha')$.
Since $\sigma$ is bijective and $\alpha'$ is a rearrangement of $\alpha$, $\sigma(\alpha')$ is a rearrangement of $\sigma(\alpha)$.
In particular, the fact that $\sigma(\alpha) \in \PF_n$ and the permutation invariance property of parking functions from \eqref{eq: characterization} together imply that $\sigma(\alpha') \in \PF_n$.
Therefore, $n + 1$ is the spot ultimately left unoccupied under $\sigma(\alpha')$, which means that $i^*$ is the spot ultimately left unoccupied under $\alpha'$.
\end{proof}

We are now ready to present our main result.
\begin{theorem}
\label{theorem: catalan}
For any $n \geq 1$, $|\PF_n^\uparrow| = C_n = \frac{1}{n + 1} \binom{2n}{n}$.
\end{theorem}
\begin{proof}
By Lemma~\ref{lemma: balls and bars} we have that $|\PP_n^\uparrow| = \binom{2n}{n}$.
Note that $\PF_n^\uparrow \subseteq \PP_n^\uparrow$.
Therefore, in the remainder of this proof we partition $\PP_n^{\uparrow}$ into equivalence classes, each of size $n + 1$ and with a single element in $\PF_n^\uparrow$.

First, let $\pi: \PP_n \rightarrow \PP_n$ be Pollak's rotation.
That is, for any $\alpha = (a_1, a_2, \ldots, a_n) \in \PP_n$, we have $\pi(\alpha) = (r_1(\alpha), r_2(\alpha), \ldots, r_n(\alpha)) \in \PP_n$ where
\begin{equation*}
    r_i(\alpha) 
    \coloneqq 
    \begin{cases}
        n + 1, & \textrm{ if } a_i + 1 \mod (n + 1) \equiv 0, \textrm{ and} \\
        a_i + 1 \mod (n + 1), & \textrm{ otherwise,}
    \end{cases}
\end{equation*}
for all $i \in [n]$.
By convention, let $\pi^{(0)} = \textrm{id}$ be the identity map and $\pi^{(s)}$ be the $s$-iterate of $\pi$.

For each parking preference tuple $\alpha \in \PP_n$, let 
\begin{equation*}
    \kappa(\alpha) = (k_1(\alpha), k_2(\alpha), \ldots, k_n(\alpha), k_{n+1}(\alpha)) \in [n]_0^{n+1}    
\end{equation*}
be its content, where
\begin{equation*}
    k_j(\alpha) \coloneqq |\{i \in [n] : a_i = j\}|
\end{equation*}
for all $j \in [n + 1]$.
In other words, $k_j(\alpha)$ is the multiplicity (or number of occurrences) of $j$ in $\alpha$.

Let $\kappa(\PP_n^{\uparrow}) = \{\kappa(\alpha) : \alpha \in \PP_n^\uparrow\} \subseteq [n]_0^{n+1}$ be the set of contents of weakly increasing parking preference tuples for the $n$ cars. 
In fact, $|\kappa(\PP_n^{\uparrow})| = |\PP_n^{\uparrow}|$ as any weakly increasing parking preference tuple is uniquely identified by its content.
Any content $k=(k_1, k_2, \ldots, k_n, k_{n+1}) \in \kappa(\PP_n^{\uparrow})$ has $n + 1$ (non-negative integer) entries and satisfies $\sum_{j = 1}^{n+1} k_i = n$, so there always exists at least one pair $j, j' \in [n+1]$ with $k_j \neq 0$ and $k_{j'} = 0$.
Therefore, the circular shift $\omega: \kappa(\PP_n^{\uparrow}) \rightarrow \kappa(\PP_n^{\uparrow})$, which moves the final entry of a tuple to the first position and shifts all other entries to the next position on the right, partitions the set $\kappa(\PP_n^{\uparrow})$ into equivalence classes, each of size $n + 1$; let $\mathcal{K}$ denote this collection of equivalence classes.
By convention, let $\omega^{(0)} = \textrm{id}$ be the identity map and $\omega^{(s)}$ be the $s$-iterate of $\omega$.

For each equivalence class $K \in \mathcal{K}$, fix an arbitrary element as its canonical representative and denote it $\hat{k} \in K$.
Moreover, for each content $k \in K$, let 
\begin{equation*}
    \tau(k) = (t_1(k), t_2(k), \ldots, t_n(k)) \in \PP_n^\uparrow    
\end{equation*}
be its corresponding weakly increasing parking preference tuple.
Note that the collection of sets $\{\tau(k): k \in K\}$ for $K \in \mathcal{K}$, each of size $n + 1$, in turn partition the set $\PP_n^\uparrow$.

Now, consider any equivalence class $K \in \mathcal{K}$, its canonical representative $\hat{k} \in K$, and its corresponding weakly increasing parking preference tuple $\tau(\hat{k}) \in \PP_n^{\uparrow}$.
For any $s \in \mathbb{N}_0$, the content of the $s$th Pollak rotation of $\tau(\hat{k})$, i.e., $\kappa(\pi^{(s)}(\tau(\hat{k})))$, is equivalently obtained as the $s$th circular shift of $\hat{k}$, i.e., $\omega^{(s)}(\hat{k})$.
Therefore, the sets $\{\pi^{(s)}(\tau(\hat{k})): s \in \mathbb{N}_0\}$ and $K$ are in bijection.
By transitivity, the sets $\{\pi^{(s)}(\tau(\hat{k})): s \in \mathbb{N}_0\}$ and $\{\tau(k) : k \in K\}$ are also in bijection, and each of them is of size $n + 1$.
Let $\phi: \{\pi^{(s)}(\tau(\hat{k})) : s \in \mathbb{N}_0\} \to \{\tau(k): k \in K\}$ be the bijective map $\phi = \tau \circ \kappa$, i.e., $\phi(\alpha) = \tau(\kappa(\alpha))$ for all $\alpha \in \{\pi^{(s)}(\tau(\hat{k})) : s \in \mathbb{N}_0\}$.
In effect, $\phi(\alpha)$ is the weakly increasing rearrangement of $\alpha$. 

Lastly, consider the sequence of $n$ cars parking on the circular one-way street with $n + 1$ numbered parking spots under any parking preference tuple $\alpha \in \{\pi^{(s)}(\tau(\hat{k})): s \in \mathbb{N}_0\} \subseteq \PP_n$ (note that, in general, $\alpha$ need not be weakly increasing).
Since there are $n + 1$ spots but only $n$ cars, all cars park and ultimately leave exactly one spot unoccupied.
By Lemma~\ref{lemma: empty spot}, this is the same spot ultimately left unoccupied if the cars arrive under the weakly increasing parking preference tuple $\phi(\alpha)$.
Moreover, exactly one preference tuple $\alpha \in \{\pi^{(s)}(\tau(\hat{k})): s \in \mathbb{N}_0\}$ results in spot $n + 1$ ultimately left unoccupied: the successive elements in $\{\pi^{(s)}(\tau(\hat{k})): s \in \mathbb{N}_0\}$ rotate the unoccupied spot clockwise around the circular street one spot at a time.
Therefore, by our previous claim, exactly one weakly increasing preference tuple $\alpha \in \{\tau(k): k \in K\}$ results in spot $n + 1$ ultimately left unoccupied.
Since the sets $\{\tau(k): k \in K\}$ for $K \in \mathcal{K}$ are each of size $n + 1$ and partition $\PP_n^\uparrow$, which is of size $\binom{2n}{n}$, there are
\begin{equation*}
    \frac{1}{n+1}\binom{2n}{n}
\end{equation*}
weakly increasing parking preference tuples that result in spot $n + 1$ ultimately left unoccupied.
These preference tuples form the set $\PF_n^\uparrow$, which establishes the desired count $|\PF_n^\uparrow| = C_n = \frac{1}{n+1}\binom{2n}{n}$.
\end{proof}

\begin{figure}[ht]
    \centering
    \includegraphics[trim={1.75cm 5cm 1.75cm 5cm},clip,width=\linewidth]{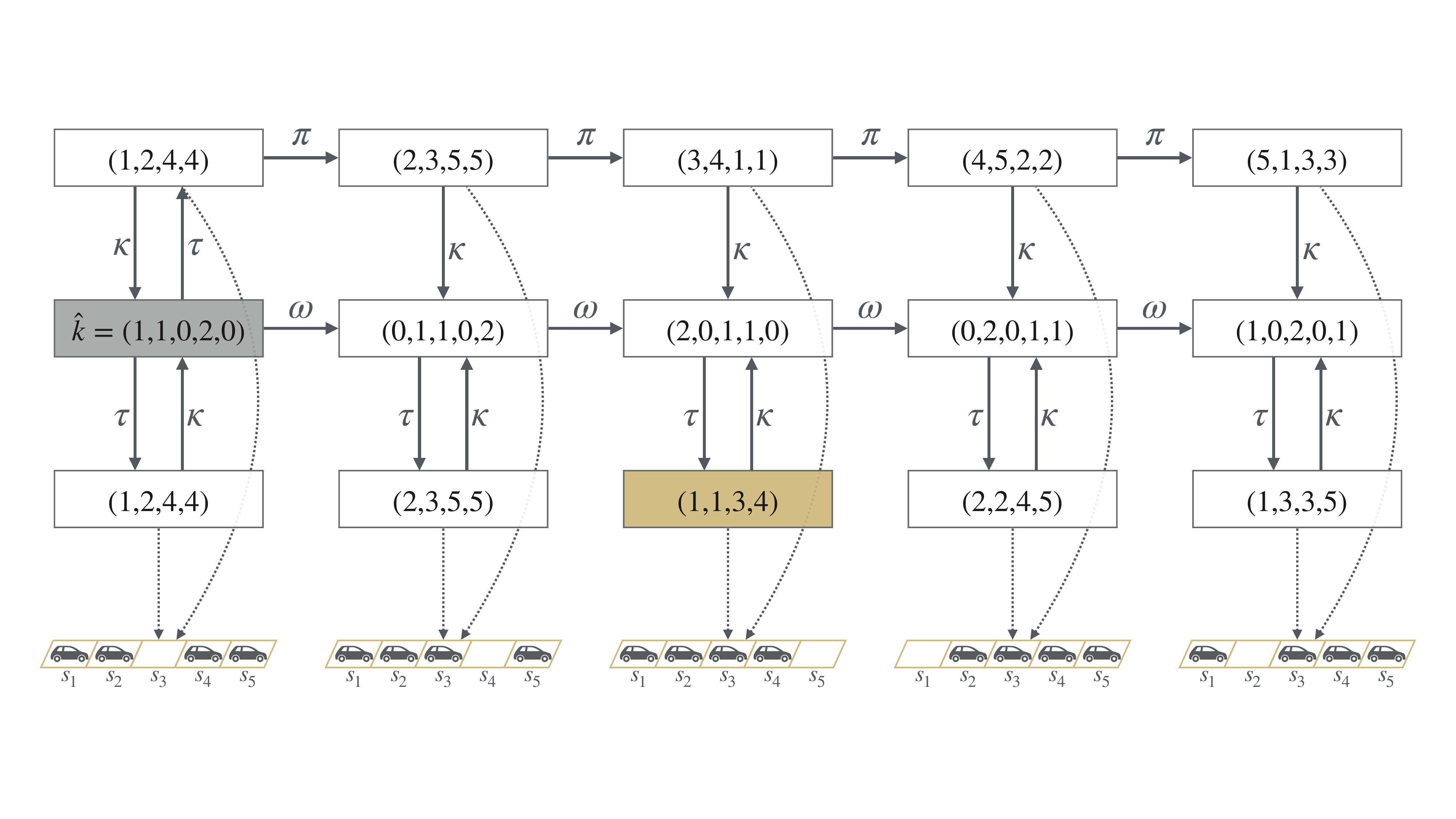}
    \caption{Illustrations accompanying \Cref{ex: proof}.
    }
    \label{fig: proof}
\end{figure}

To summarize: 
\begin{enumerate}
    \item 
    Lemma~\ref{lemma: balls and bars} establishes the $\binom{2n}{n}$ count for the broader set of weakly increasing parking preferences.
    \item 
    Given this broader set, we would like to use Pollak's rotation $\pi$ to partition it into equivalence classes, each of size $n + 1$ and with a single weakly increasing parking function.
    However, $\pi$ need not preserve the weakly increasing property, so it is not immediately clear how to form these.
    To do so, we start with any weakly increasing parking preference and:
    \begin{itemize}
        \item 
        Send its Pollak rotations to their weakly increasing rearrangements via the map $\phi$.
        The content map $\kappa$ and the weakly increasing tuple map $\tau$ make this precise: $\kappa$ links the Pollak rotations to circular shifts $\omega$ of the initial content, while $\tau$ in turn links these contents to the weakly increasing rearrangements.
        That is, $\phi = \tau \circ \kappa$.
        \item 
    	Among these Pollak rotations, there is a single parking function (not necessarily weakly increasing).
    	Lemma~\ref{lemma: empty spot} implies that the weakly increasing rearrangement of this unique parking function is the unique weakly increasing parking function in its equivalence class. 
    \end{itemize}
    Note that in the enumeration of parking functions, each equivalence class is the set of all Pollak rotations of a given preference tuple.
    Conversely, in our enumeration of weakly increasing parking functions, each equivalence class is instead formed by circular shifts of the content of a given weakly increasing preference tuple: its Pollak rotations are instead used to keep track of the empty spot.
    For example, as seen in Example~\ref{ex: proof} and its accompanying Figure~\ref{fig: proof}, (2,2,4,5) is not a Pollak rotation of (1,2,4,4), yet (1,2,4,4) and (2,2,4,5) are in the same equivalence class.
\end{enumerate}

In \Cref{ex: proof} and its accompanying \Cref{fig: proof} we illustrate the main ideas in our proof of \Cref{theorem: catalan}.
\begin{example}    
\label{ex: proof}
    Consider $4$ cars parking on a circular one-way street with $5$ spots, i.e., after spot $5$, the street loops back into spot $1$.
    The weakly increasing preferences of length $4$ are in bijection with the contents of length $5$ via the content map $\kappa$ whose inverse is the weakly increasing word map $\tau$.
    For example, $\kappa\left((1,1,3,4)\right) = (2,0,1,1,0)$ and $\tau\left((2,0,1,1,0)\right) = (1,1,3,4)$.
    
    The circular shift operation on tuples, denoted by $\omega$, partitions the set of contents into equivalence classes of size $5$.
    In particular, the elements of the orbit
    \begin{equation*}
        (1,1,0,2,0) \stackrel{\omega}{\rightarrow} 
        (0,1,1,0,2) \stackrel{\omega}{\rightarrow} 
        (2,0,1,1,0) \stackrel{\omega}{\rightarrow} 
        (0,2,0,1,1) \stackrel{\omega}{\rightarrow} 
        (1,0,2,0,1) \stackrel{\omega}{\rightarrow} 
        (1,1,0,2,0)
    \end{equation*}
    are all in the same equivalence class: without loss of generality, let $\hat{k} = (1,1,0,2,0)$ be its canonical representative (shaded in grey in \Cref{fig: proof}).
    We also rotate its corresponding weakly increasing preference $\tau(\hat{k}) = (1,2,4,4)$ via the rotation map $\pi$, which increases all of its entries by one modulo $n + 1$, as in Pollak's original proof.
    This forms the orbit
    \begin{equation*}
        (1,2,4,4) \stackrel{\pi}{\rightarrow} 
        (2,3,5,5) \stackrel{\pi}{\rightarrow} 
        (3,4,1,1) \stackrel{\pi}{\rightarrow} 
        (4,5,2,2) \stackrel{\pi}{\rightarrow}
        (5,1,3,3) \stackrel{\pi}{\rightarrow} 
        (1,2,4,4),
    \end{equation*}
    also of size $5$, and with exactly one (not necessarily weakly increasing) parking function: $(3,4,1,1)$.
    Through $\kappa$, these preference rotations are in bijection with the equivalence class of contents.
    For example, 
    \begin{equation*}
        \kappa\left((3,4,1,1)\right) = (2,0,1,1,0).
    \end{equation*}
    Moreover, the composition $\phi = \tau \circ \kappa$ bijectively maps these preference rotations to weakly increasing preferences.
    For example, 
    \begin{equation*}
        \phi\left((3,4,1,1)\right) = \tau\left((2,0,1,1,0)\right) = (1,1,3,4).
    \end{equation*}
    Finally, the spot ultimately left unoccupied is invariant under $\phi$, as noted by the dashed arrows and the unlabeled parking outcomes.
    For example, the spot ultimately left unoccupied under both $(3,4,1,1)$ and $(1,1,3,4)$ is spot $5$.
    Since $(3,4,1,1)$ is the unique parking function in its equivalence class of preference rotations, $(1,1,3,4)$ is the unique parking function in its equivalence class of weakly increasing preferences (shaded in golden in \Cref{fig: proof}).

\end{example}
\section{Discussion}
\label{sec: discussion}

Circular-street arguments have been useful in combinatorics specific to parking-related objects.
\citet{durmic2023probabilistic} used a circular-street argument to show that, under a coin-flipping parking protocol in which each car moves forward or backward depending on its coin-flipping outcome, the probability that a parking preference tuple sampled uniformly at random is a parking function is independent of the coin bias parameter $0 \leq p \leq 1$.
More recently, \citet{rubey2025fixed} developed circular-street arguments to enumerate certain permutation-related discrete statistics of parking functions.
In a similar spirit, our proof of Theorem~\ref{theorem: catalan} further demonstrates the versatility of circular-street arguments as a general proof technique.

\acknowledgements
\label{sec:ack}
We thank the anonymous referees for suggestions that improved our level of exposition.

\nocite{*}
\bibliographystyle{abbrvnat}
\bibliography{circular}

\begin{thebibliography}{13}
\providecommand{\natexlab}[1]{#1}
\providecommand{\url}[1]{\texttt{#1}}
\expandafter\ifx\csname urlstyle\endcsname\relax
  \providecommand{\doi}[1]{doi: #1}\else
  \providecommand{\doi}{doi: \begingroup \urlstyle{rm}\Url}\fi

\bibitem[Armstrong et~al.(2016)Armstrong, Loehr, and Warrington]{armstrong2016rational}
D.~Armstrong, N.~A. Loehr, and G.~S. Warrington.
\newblock Rational parking functions and {C}atalan numbers.
\newblock \emph{Ann. Comb.}, 20\penalty0 (1):\penalty0 21--58, 2016.
\newblock ISSN 0218-0006,0219-3094.
\newblock URL \url{https://doi.org/10.1007/s00026-015-0293-6}.

\bibitem[Durmi\'c et~al.(2023)Durmi\'c, Han, Harris, Ribeiro, and Yin]{durmic2023probabilistic}
I.~Durmi\'c, A.~Han, P.~E. Harris, R.~Ribeiro, and M.~Yin.
\newblock Probabilistic parking functions.
\newblock \emph{Electron. J. Combin.}, 30\penalty0 (3):\penalty0 Paper No. 3.18, 25, 2023.
\newblock ISSN 1077-8926.
\newblock URL \url{https://doi.org/10.37236/11649}.

\bibitem[Dvoretzky and Motzkin(1947)]{dvoretzky1947problem}
A.~Dvoretzky and T.~Motzkin.
\newblock A problem of arrangements.
\newblock \emph{Duke Math. J.}, 14:\penalty0 305--313, 1947.
\newblock ISSN 0012-7094,1547-7398.
\newblock URL \url{http://projecteuclid.org/euclid.dmj/1077474127}.

\bibitem[Golomb(1956)]{golomb56}
S.~W. Golomb.
\newblock Classroom {N}otes: {C}ombinatorial {P}roof of {F}ermat's "{L}ittle" {T}heorem.
\newblock \emph{Amer. Math. Monthly}, 63\penalty0 (10):\penalty0 718, 1956.
\newblock ISSN 0002-9890,1930-0972.
\newblock URL \url{https://doi.org/10.2307/2309563}.

\bibitem[Konheim and Weiss(1966)]{konheim1966occupancy}
A.~G. Konheim and B.~Weiss.
\newblock An occupancy discipline and applications.
\newblock \emph{SIAM Journal on Applied Mathematics}, 14\penalty0 (6):\penalty0 1266--1274, 1966.
\newblock \doi{10.1137/0114101}.

\bibitem[Mart\'inez~Mori(2024)]{martinezmori2025what}
J.~C. Mart\'inez~Mori.
\newblock What is{$\ldots$} a parking function?
\newblock \emph{Notices Amer. Math. Soc.}, 71\penalty0 (8):\penalty0 1062--1065, 2024.
\newblock ISSN 0002-9920,1088-9477.

\bibitem[{OEIS Foundation Inc.}()]{OEIS}
{OEIS Foundation Inc.}
\newblock The on-line encyclopedia of integer sequences.
\newblock \url{https://oeis.org}.
\newblock Accessed: 2026-04-08.

\bibitem[Pyke(1959)]{pyke1959supremum}
R.~Pyke.
\newblock The supremum and infimum of the poisson process.
\newblock \emph{The Annals of Mathematical Statistics}, 30\penalty0 (2):\penalty0 568--576, 1959.
\newblock \doi{10.1214/aoms/1177706269}.

\bibitem[Raney(1960)]{raney1960functional}
G.~N. Raney.
\newblock Functional composition patterns and power series reversion.
\newblock \emph{Trans. Amer. Math. Soc.}, 94:\penalty0 441--451, 1960.
\newblock ISSN 0002-9947,1088-6850.
\newblock URL \url{https://doi.org/10.2307/1993433}.

\bibitem[Riordan(1969)]{riordan1969ballots}
J.~Riordan.
\newblock Ballots and trees.
\newblock \emph{Journal of Combinatorial Theory}, 6\penalty0 (4):\penalty0 408--411, 1969.
\newblock \doi{10.1016/S0021-9800(69)80039-6}.

\bibitem[Rubey and Yin(2025)]{rubey2025fixed}
M.~Rubey and M.~Yin.
\newblock Fixed points and cycles of parking functions.
\newblock \emph{Enumer. Combin. Appl.}, 5\penalty0 (2):\penalty0 Article \#S2R10, 2025.

\bibitem[Schumacher(2018)]{schumacher2018descents}
P.~R.~F. Schumacher.
\newblock Descents in parking functions.
\newblock \emph{J. Integer Seq.}, 21\penalty0 (2):\penalty0 Art. 18.2.3, 8, 2018.
\newblock ISSN 1530-7638.

\bibitem[Stanley(2015)]{stanley2015catalan}
R.~P. Stanley.
\newblock \emph{Catalan Numbers}.
\newblock Cambridge University Press, 2015.
\newblock \doi{10.1017/CBO9781139871495}.

\end{thebibliography}
\label{sec:biblio}

\end{document}